\documentclass[12pt]{article}

\usepackage{amsfonts}
\usepackage{amssymb}
\usepackage{amsmath}
\usepackage{amsthm}
\usepackage{graphicx}
\usepackage{float}
\usepackage{bm}
\usepackage{color}
\usepackage{subcaption}
\usepackage{cite}

\newcommand{\bsx}{{\bf x}}

\newtheorem{exmp}{Example}[section]

\theoremstyle{plain}

\theoremstyle{remark}
\newtheorem{remark}{Remark}

 \title{Fourier smoothed pre-corrected trapezoidal rule \\ for solution of Lippmann-Schwinger integral equation}
\author{Ambuj Pandey \footnote{Indian Institute of Science Education and Research,
		Bhopal, India, Email: ambuj@iiserb.ac.in} \and Akash Anand\footnote{Indian Institute of Technology Kanpur, Kanpur, India, Email: akasha@iitk.ac.in}} \date{}

\begin{document}
	\maketitle
		

\section{Introduction}\label{sec:Introduction}
The numerical simulation of wave scattering by penetrable inhomogeneous media is of paramount importance in many fields of science and engineering, for instance, medical imaging, underwater acoustics, radar, remote sensing, and plasma physics, to name a few. Formally, in the two dimensions, the problem we consider in this paper reads as follows: for a given incident field $u^{i}$ with wavenumber $\kappa$ that satisfy free space Helmholtz equation $\Delta u^i+\kappa^2u^i=0$ in $\mathbb{R}^2$ and a bounded inhomogeneity $D \subset\mathbb{R}^2$ with a variable refractive index $\mu(\bsx)$, find the scattered field $u^s$ satisfying the Sommerfeld radiation condition
\begin{equation}\label{SRC}
\lim_{r\to\infty} r^{1/2}\left(\frac{\partial u^s}{\partial r}-i\kappa u^s\right)=0, 
\end{equation}
where $ r=|{\bf x}|$ is the Euclidean length of ${\bf x} = (x_1,x_2) \in \mathbb{R}^2$, such that the total field  $u = u^{s} + u^{i}$ satisfies
 \begin{equation}\label{utotal}
 \Delta u+\kappa^2\mu^2u=0, \qquad \text{in\ } \mathbb{R}^2.
 \end{equation}
In the homogeneous region outside the obstacle $D$, the refractive index $\mu \equiv 1$ and, in general, is discontinuous across the boundary $\partial D$. 
An equivalent integral equation formulation corresponding to the problem~\eqref{SRC}-\eqref{utotal}  is given by Lippmann-Schwinger equation~\cite{martin2003acoustic},
\begin{equation}\label{Lipp}
u({\bf x})+\kappa^2\int_{\mathbb{R}^2}G({\bf x-y})m({\bf y})u({\bf y})d{\bf y}=u^i({\bf x}), \qquad {\bf x}\in\mathbb{R}^2,  
\end{equation}
with the contrast function $m(\bsx)=1-\mu^2(\bsx)$ and  $G({\bf x})=i H_0(\kappa|{\bf x}|)/4$ where  we  denote the Hankel function of the first kind of order zero by $H_0$. 

One of the main advantages of using equation~\eqref{Lipp} is that radiation condition~\eqref{SRC} is automatically satisfied. Moreover, unlike the schemes based on the differential formulation~\eqref{SRC}-\eqref{utotal} or their variational counterparts, a typical Nystr\"{o}m discretization of~\eqref{Lipp} produces a linear system for which the condition number does not grow with grid refinements and the corresponding numerical solution does not suffer from the dispersion error.
However, obtaining such a discretization that is computationally efficient as well as numerically accurate presents a couple of tough design challenges.  The first significant difficulty pertains to the accurate numerical integration of the convolution operator in \eqref{Lipp} whose integrand is singular within the integration domain owing to the factor $G({\bf x-y})$.    
Achieving a reduced computational complexity poses the other challenge. We see that the resulting dense and non-Hermitian linear system forces the use of matrix-free iterative schemes like GMRES for its solution. Thus, a straightforward integration scheme requiring $O(N^2)$ operations at each iteration can get prohibitively expensive with increasing $\kappa \mathrm{d} \Vert \mu\Vert_{\infty}$ as the discretization size $N$ would need to grow as $O((\kappa  \mathrm{d} \Vert \mu \Vert_{\infty})^2)$ to maintain a fixed accuracy for inhomogeneity of diameter $ \mathrm{d} $. The goal, therefore, is to develop an accurate integration scheme for the approximation of the convolution operator in equation~\eqref{Lipp} with computational complexity $O(N \log N)$.

Several numerical schemes in this direction are available in the literature.  While we do not intend to review all such work here, some relevant contributions can be found in~\cite{anand2016efficient,pandey2019improved,aguilar2004high,bruno2004efficient,
	duan2009high,hyde2002fast,andersson2005fast,zepeda2016fast, 	liu2018sparsify,sifuentes2018preconditioning} and references therein.  Among all existing fast methods, pre-corrected quadrature rule~\cite{aguilar2004high,duan2009high} is perhaps the simplest in terms of implementation. This scheme basically computes, using FFT, a two dimensional Toeplitz convolution on an underlying uniform rectangular grid. While method converges to high-order for smooth compactly supported densities, it exhibits only the first order convergence in the case when the contrast function $m$ is discontinuous across the material interface $\partial D$ that we consider in this paper. In this short article, we utilize a ``Fourier smoothing technique" which, when used in conjunction with the pre-corrected quadrature rule~\cite{aguilar2004high,duan2009high}, yields a second order convergence for discontinuous $m$ while maintaining its original computational complexity of $O(N \log N)$. The resulting algorithm, which we refer to as ``Fourier smoothed pre-corrected trapezoidal rule", is still very simple to implement and its performance is shown to exceed that of other well-known $O(N \log N)$ methods that converge quadratically for a discontinuous material interface. Moreover, the method is not limited to simple or smooth scattering configurations and converges with second order even for geometrically complex discontinuous scatterers, such as inhomogeneities with corners and cusps, without adding any additional component in the overall algorithm. The numerical results in Section~\ref{CR} clearly demonstrate that this method performs better in comparison with some of its popular alternatives. 
The rest of this short note is organized as follows: we present the Fourier smoothed pre-corrected trapezoidal rule in Section~\ref{WFSM}.  To corroborate the second order convergence and to make other relevant performance comparisons, we include a variety of numerical experiments in Section~\ref{CR}. Finally, we make some concluding remarks in Section~\ref{CL}.

\section{Fourier smoothed pre-corrected trapezoidal rule}\label{WFSM}
\begin{figure}[t!] 
  	\begin{center}	
  			\begin{minipage}[b]{.46\linewidth}
  			\centering
  			\includegraphics[width=\linewidth]{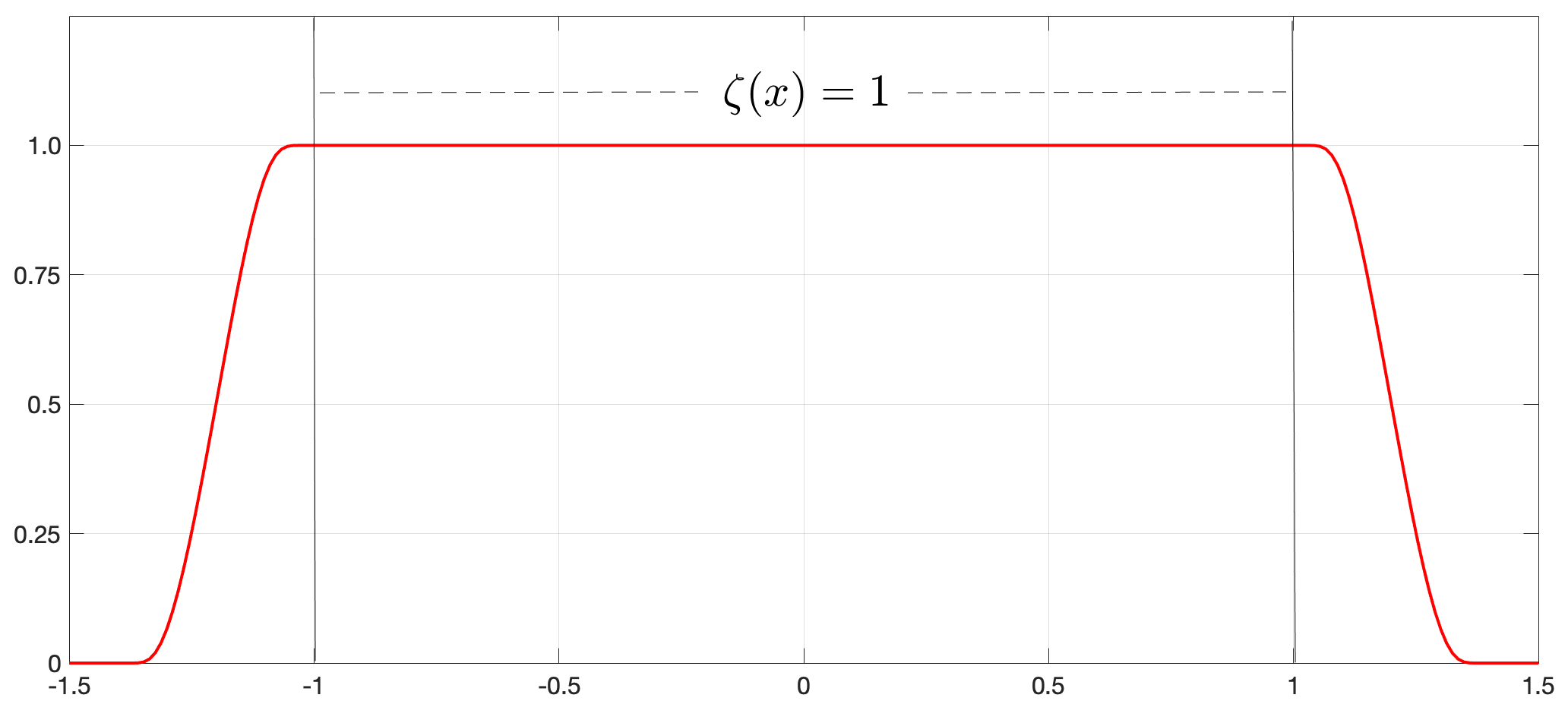}
  			\subcaption{One dimensional smooth vanshing window function $\zeta(x)$ with $r_{\mathrm{in}}=1.0$ and $r_{\mathrm{out}}=1.4$ . Function  $\zeta(x)$  and all of its derivative 
  			vanishes along the end points of interval $[-1.5,1.5]$.  }
  			\label{fig:1d_window}
  		\end{minipage}
  	\hspace{7mm}
  		\begin{minipage}[b]{.48\linewidth}
  			\centering
  			\includegraphics[width=\linewidth]{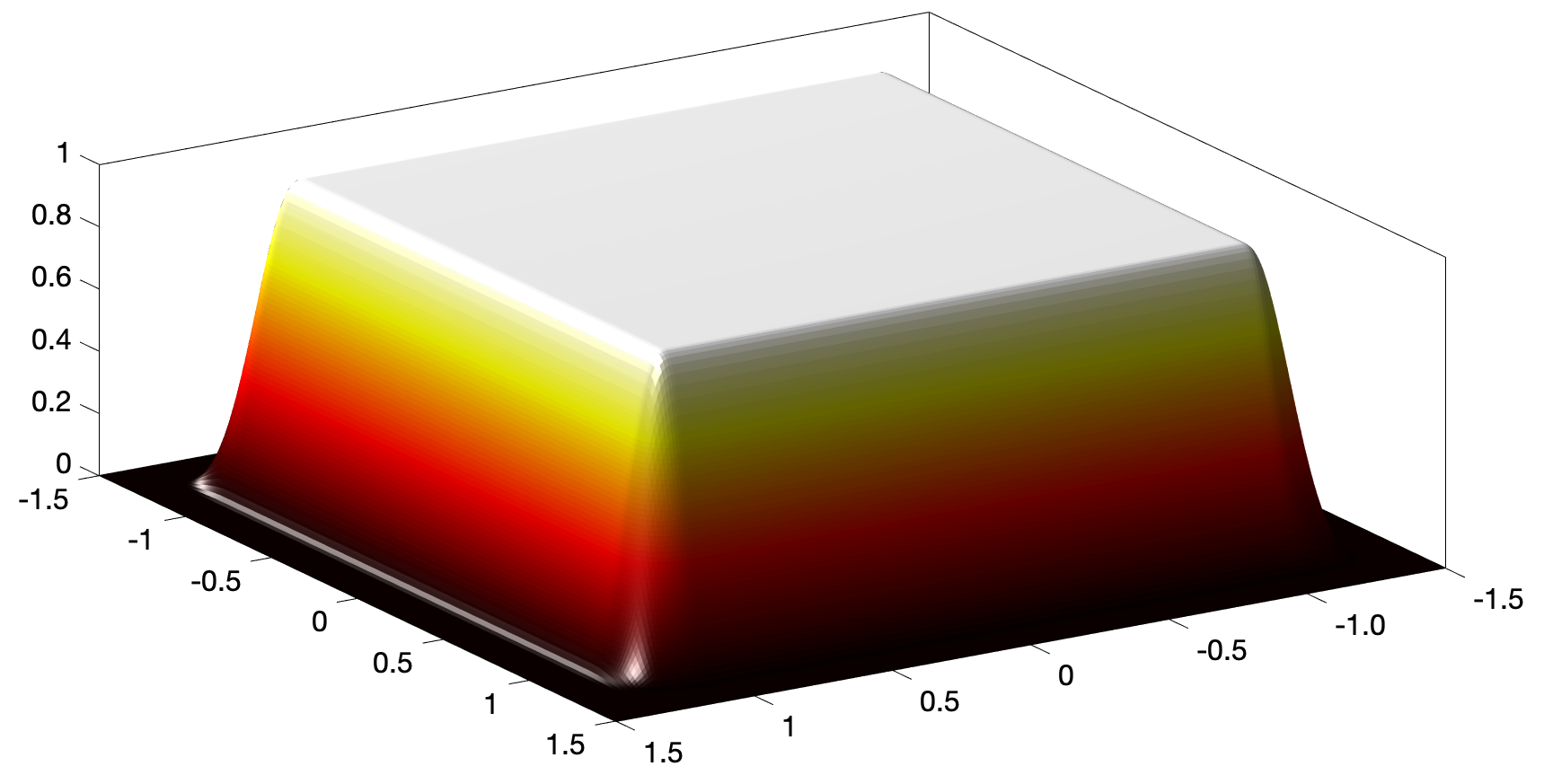}
  			\subcaption{Two dimensional smooth vanshing window function $\eta$ which is obtained as $\eta(\bsx)=\zeta(x_1)\zeta(x_2)$ and equal to one if $\bsx \in [-1,1]^2$ and vanishes to zero towrds the boundary of $[1.5,1.5]^2$.  }
  			\label{fig:2d_window}
  		\end{minipage}  		
  	\end{center}	
  \end{figure}

We present  the details in the context of
 approximation of convolution operator $\mathcal{K}$ given by
 \begin{equation}\label{CO}
\mathcal{K}(u)(\bsx) =\int_{\Omega}G({\bf x-y})m({\bf y})u({\bf y})d{\bf y}, \ \ \ \bsx \in \Omega,
\end{equation}  
where the square $\Omega=[-a,a]^{2}$, $a>0$, is such that  $\Omega$ contains the obstacle $D$.
Here, motivated by the integral operator in the Lippmann-Schwinger integral equation \eqref{Lipp}, we assume that $u$ is known on $\Omega$ while $m \equiv 0$ outside $D$. We also assume that $m$ on $D$ has a smooth extension on $\Omega$, that is, there exits a smooth function $\tilde{m}$ defined on $\Omega$ such that $\tilde{m} \equiv m$ on $D$. This allows for construction of a bi-periodic Fourier extension of first kind $m_e$ of $m$ \cite{boyd} by simply multiplying $\tilde{m}$ by a smooth window function $\eta$ such that $\eta({\bf x})=1$ if ${\bf x} \in D$ and vanishes to high-order at the boundary of $\Omega$.  Note that such a window function $\eta$ is not difficult to construct and can be easily obtained as a tensor product of one dimensional smooth window functions. For instance, in all numerical results presented in Section~\ref{CR}, we have used $\eta(\bsx) =\zeta(x_1)\zeta(x_2)$ where 
\begin{equation}\label{SWF}
\zeta(x)=\begin{cases}
1, & \text{for } |x|\le r_{\mathrm{in}}\\
\exp{(2e^{-1/r}/(r-1))}, & \text{for } r_\mathrm{in}<|x|<r_\mathrm{out}, \text{where }r=\frac{|x|-r_\mathrm{in}}{r_\mathrm{out}-r_\mathrm{in}}\\
0, & \text{for } |x|\ge r_{\mathrm{out}},
\end{cases}
\end{equation}	
for appropriately chosen $r_{\mathrm{in}}$ 
and $r_{\mathrm{out}}$. 
An example $\zeta$ and $\eta $ are shown in Figures~\ref{fig:1d_window} and~\ref{fig:2d_window} respectively.

Now, using an indicator function $\chi_{D}$ that takes value one in $D$ and zero otherwise, we recast the integral operator~\eqref{CO} as
\begin{equation} \label{WLS}
 \mathcal{K}(u)(\bsx) =\int_{\Omega}G({\bf x-y})m({\bf y})u({\bf y})\chi_{D}({\bf y})d{\bf y} =\int_{\Omega}G({\bf x-y})m_e({\bf y})u({\bf y})\chi_{D}({\bf y})d{\bf y}.
\end{equation}
Motivated by the Fourier smoothing used in~\cite{hyde2005fast} to soften the discontinuity in the contrast function, we work with a truncated Fourier series $ \chi_{D}^{F}({\bf x})$ of $ \chi_{D}({\bf x})$ of period $2a$ of the form
\begin{equation}\label{FCI}
\chi_{D}^{F}({\bf x}) =\sum_{j=-F}^{F}\sum_{k=-F}^{F}c_{jk} \exp\left({\frac{\pi ijx_1}{a} }\right)\exp\left({\frac{\pi ikx_2}{a} }\right),
\end{equation}
where the Fourier coefficients $c_{jk}$ are assumed to be known analytically or can be computed very accurately.  The simplicity of this procedure lies in the fact that, unlike in~\cite{hyde2005fast},  the smoothing is applied independent of the density in~\eqref{CO} to a piecewise constant function that depends only on $D$. The proposed Fourier smoothed solver, therefore, can seamlessly handle any variable refractive index $\mu(\bsx)$ on a given $D$ once the Fourier coefficients in~\eqref{FCI} are obtained for that particular geometry.
 

   
Our procedure for approximation of~\eqref{WLS} at any grid point ${\bf x}_{j}\in \Omega_{n}$, where
\[
\Omega_{n} =\left\{ \bsx_{j}=-(a,a)+h(j_1,j_{2})\ |\ 0\leq j_1, j_{2} < n, h=2a/n \right\},
\]  
is a uniform computational grid on $\Omega$ of size $N = n^2$,  
is composed of two main steps: first, we Fourier smooth the discontinuity in $\chi_D$ in~\eqref{WLS} by replacing it with $\chi_{D}^{F}$ , and second, approximate the resulting integral using the pre-corrected trapezoidal rule leading to the approximation
\begin{equation}\label{CTR}
 \mathcal{K}(u)(\bsx_j) \approx \int_{\Omega}G({{\bf x}_{j},{\bf y}})m_e({\bf y})u({\bf y})\chi_{D}^{F}({\bf y})d{\bf y}\approx \sum_{{\bf x}_{k} \in {\Omega_{n}}}w_{|j-k|}m_e({\bf x}_k)u({\bf x}_k)\chi_{D}^{F}({\bf x}_k).
\end{equation}   
Here, weights $w_{| j-k |}$ are given by
\begin{equation}\label{DC_2d}
w_{| j-k |}=\begin{cases}
h^{2} G({\bsx_{j}-\bsx_k}),& \text{if } j \neq k, \\
\frac{h^2}{2\pi}\left[\pi i/2-\left(\ln\left(\kappa h/2\right) -\gamma-\beta \right)\right], & \text{otherwise},
\end{cases}
\end{equation} 
where 
$\gamma$ 
is the Euler's constant and $\beta$ is the corrected coefficient for Green's function singularity whose numerical value can be read from the Table 1 in \cite{aguilar2004high} or can be computed accurately using the analytical expressions given in~\cite{duan2009high}.  It is important to note that, for the quadratic convergence as $N(=n^2)$ increases, we also need to proportionally increase the truncation size $F$ in~\eqref{FCI}. For instance, we use $F=n/2$ for the computations presented in this text. The evaluation of $\chi^{F}_{D}$ at all grid points  $\bsx_j \in \Omega_{n}$ can be done in $O(N \log N)$ operations by means of FFT. Further, as the discrete summation on the right-hand side of~\eqref{CTR} is a two dimensional Toeplitz convolution, it can also be computed at all grid points $\bsx_j \in \Omega_{n}$  in $O(N \log N)$ operations by means of FFT. The overall computational cost of the algorithm, thus, remains at $O(N \log N)$.
\begin{remark}
While we presented the Fourier smoothing of pre-corrected trapezoidal rule, as evident from \eqref{CTR} where we see that application of Fourier smoothing and singular quadratures are independent steps, we expect that it can similarly be
combined with other existing linearly convergent fast solvers when dealing with discontinuous scatterers to enhance their speed of convergence.
\end{remark}
\begin{remark}
This approach can obviously be adapted to other problems where the convolution integral of the form~\eqref{CO} with discontinuous density is required to be computed. In such a scenario, depending on the underlying singularity of the kernel $G$, the corresponding pre-corrected quadrature would be required. In the case when the kernel is smooth, the classical trapezoidal rule replaces the pre-corrected quadrature.
\end{remark}
%
\section{Computational Results} \label{CR}
In this section, we present numerical evidence for the quadratic convergence of the proposed Fourier smoothed pre-corrected trapezoidal (FSPT) rule through a variety of numerical experiments involving discontinuous inhomogeneities including those with corners, and cusps. The reported errors are measured using the expressions
\begin{align*}
\varepsilon_{\infty}= \frac{ \underset{1\leq j \leq N}{\max} \left|u^{\text{exact}}({\bsx}_{j})-u^{\text{approx}}(\bsx_{j})\right|}{ \underset{1\leq j \leq N}{\max} \left|u^{\text{exact}}(\bsx_{j}) \right|}, \ \  \text{and}\ \
&	\varepsilon_{2}= \left ( \frac{\sum \limits_{j=1}^{N}\left|u^{\text{exact}}(\bsx_{j})-u^{\text{approx}}(\bsx_{j})\right|^2}{\sum \limits_{j=1}^{N}\left|u^{\text{exact}}(\bsx_{j}) \right|^2} \right)^\frac{1}{2},
\end{align*}    
where $\bsx_{j} \in \Omega_n$. In those cases where the analytical solution is not available, $u^{\text{exact}}$ in above expressions is taken to be the numerical solution obtained by FSPT for a computational grid finer than those included in that particular convergence study. 
The symbol ``F" in tables reported in this section represents the number of Fourier modes in equation~\eqref{FCI} used for the corresponding Fourier smoothing while ``noc'' denotes the numerical order of convergence. Moreover,
in all examples discussed in this section, the contrast function $m$ is taken to be discontinuous across the interface $\partial D$ with $m \equiv 0$ outside $D$. 

\begin{exmp}(\textit{Scattering by a circular inhomogeneity})	
\end{exmp}
As a first example, we take the canonical problem of scattering by a circular scatterer  $D$ of unit radius. For this experiment, the contrast function is taken to be  $m(\bsx)=-0.5$ for $\bsx\in D$. The computational domain used is $\Omega=[-1.1,1.1]\times[-1.1,1.1]$ while we use $r_{\mathrm{in}}=1.01$ and $r_{\mathrm{out}}=1.08$ in the construction of $\eta$.  The incoming wave is taken as $u^{i}(\bsx)=J_{0}(\kappa |{\bsx}|)$ with $\kappa = 10$. The analytical solution, in this case, is available (for example, see~\cite{andersson2005fast}). The corresponding errors are reported in Table~\ref{tabel:-Disc}. These results clearly demonstrate the significant improvement in accuracy and the enhancement in the rate of convergence for the proposed Fourier smoothed pre-corrected method over the direct pre-corrected scheme. For instance, the last row of Tables~\ref{tabel:-Disc} indicates a thousand times better accuracy in comparison.

\begin{table}[t!]
	\begin{center} 
		\begin{tabular}{c|c|c|c|c|c|c|c|c|c} \hline
			\hline		
			$n$ & \multicolumn{4}{c|}{pre-corrected} &\multicolumn{5}{c}{Fourier smoothed pre-corrected}  \\ 
			\cline{2-10} &$\varepsilon_{2} $& noc& $\varepsilon_{\infty}$& noc &F &  $\varepsilon_{2}$& noc & $\varepsilon_{\infty}$& noc \\ 
			\hline
			\hline
			$16$&$9.7\times 10^{-2}$& - &$8.8\times 10^{-2}$ & -&8 & $7.5 \times 10^{-2}$&-&$5.2\times 10^{-2}$\\ 
			
			$32$ &$6.7\times 10^{-3}$ & 3.7&$4.6\times 10^{-3}$ &4.2 & 16 &$3.9\times 10^{-3}$&4.3& $2.8 \times 10^{-3}$& 4.2\\ 
			$64$ &$6.0\times 10^{-3}$ &0.2&$3.9\times 10^{-3}$&0.3 &32 &$1.3\times 10^{-3}$&1.6&$7.6 \times 10^{-4}$&1.9\\ 
			$128$&$1.9\times 10^{-3}$& 1.6 &$1.3\times 10^{-3}$&1.6 & 64 &$2.7\times 10^{-5}$&5.6&$3.2\times 10^{-5}$& 4.6\\ 
			$256$ &$7.3\times 10^{-4}$ & 1.4& $5.2\times 10^{-4}$& 1.3&128&$3.7\times 10^{-6}$&2.9&$5.8 \times 10^{-6}$&2.5\\ 
			
			$512$ &$2.5\times 10^{-4}$& 1.5&$1.3\times 10^{-4}$ &1.9 &256&$5.4\times 10^{-7}$&2.8&$1.2 \times 10^{-6}$& 2.3\\ 
			
			$1024$& $8.5\times 10^{-5}$& 1.5&$5.7\times 10^{-5}$ & 1.3&512 &$9.0\times 10^{-8}$&2.6&$2.7 \times 10^{-7}$& 2.1\\ 
			
			$2048$& $3.9\times 10^{-5}$& 1.1 & $3.8\times 10^{-5}$& 0.6& 1024&$1.6\times 10^{-8}$&2.5&$6.4 \times 10^{-8}$& 2.1\\ 				
			\hline	
			\hline
		\end{tabular}  
		\caption{A convergence report for the solution of Lippmann-Schwinger equation~\eqref{Lipp} with and without Fourier smoothing in the $\varepsilon_2$ measure. The incidence wave used is $u^{i}(\bsx)=J_{0}(10|\bsx|)$ and the contrast function $m$ takes the value $-0.5$ inside the circular scatterer centered at the origin. The errors are obtained by comparing the numerical solutions against the known analytical solution.} 
		\label{tabel:-Disc} 
	\end{center} 
\end{table}

\begin{table}[b!]
	\begin{center} 
		\begin{tabular}{c|c|c|c|c|c|c} \hline
			\hline
			Method & Unknowns & F &$\kappa$&\multicolumn{2}{c|}{Error}& {Time}\\ 
			\cline{5-6}& &  & &  $\varepsilon_{\infty}$   & $\varepsilon_{2}$& (in sec.)   \\ 
			\hline
			\hline

			AH &$65536$& - &$100\pi$ &$4.4\times 10^{-3}$&$4.0\times 10^{-4}$&53.1\\ 
			AH &$65536$& -&$120\pi$ &$6.0\times 10^{-2}$&$8.3\times 10^{-3}$&77.1\\ 

			FSPT &$65536$& 128 &100$\pi$ &$2.5\times 10^{-3}$&$3.3\times 10^{-3}$&1.0\\ 	
			FSPT &$65536$& 128 &120$\pi$ &$9.3\times 10^{-3}$&$1.2\times 10^{-3}$&2.0\\ 
			
			FSPT &$1048576$& 512 &120$\pi$ &$1.1\times 10^{-5}$&$8.4\times 10^{-6}$&44.0\\

			\hline	
			\hline
		\end{tabular}  
		\caption{Comparison of Fourier smoothed pre-corrected trapezoidal (FSPT) with the solver by Anderson and Holst (AH) reported in \cite{andersson2005fast}.} 
		\label{AH} 
	\end{center} 
\end{table} 

\begin{table}[t!]
	\begin{center} 
		\begin{tabular}{c|c|c|c|c|c|c} \hline
			\hline
			Method & Unknowns & F &$\kappa$&\multicolumn{2}{c|}{Error}& {Time}\\ 
			\cline{5-6}& &  & &  $\varepsilon_{\infty}$   & $\varepsilon_{2}$& (in sec.)  \\ 
			\hline
			\hline
			
			BH &$262000$& - &$2\pi$ &$5.6\times 10^{-5}$&-&99\\ 
			BH  &1047000& - &$2\pi$&$2.0\times 10^{-6}$&-&808\\ 
			
			FSPT &$262144$& 256 &$2\pi$ &$9.0\times 10^{-6}$&$1.0\times 10^{-6}$ &3.0\\ 	
			FSPT &$1048576$& 512 &2$\pi$ &$1.9\times 10^{-6}$&$1.9\times 10^{-7}$&15.0\\ 		
			\hline	
			\hline
		\end{tabular}  
		\caption{Comparison of Fourier smoothed pre-corrected trapezoidal (FSPT) with the method by Bruno and Hyde (BH) reported in \cite{hyde2002fast}.} 
		\label{BH} 
	\end{center} 
\end{table}

\begin{exmp}(\textit{Comparison with other second order methods})
\end{exmp}
This example compares the accuracy and run time performance of the Fourier smoothed pre-corrected trapezoidal scheme with a couple of other existing fast, second order convergent methods. In particular, we compare with the solver reported in~\cite{andersson2005fast} by Anderson and Holst (AH) and with the approach developed by Bruno and Hyde (BH) in~\cite{bruno2004efficient, hyde2002fast}.  

A comparative study with the AH scheme is presented in Table~\ref{AH}. In this case, the solution of  equation~\eqref{Lipp} is computed for the incidence wave $u^{i}(\bsx)=J_{0}(\kappa |\bsx|)$, and contrast function $m(\bsx)=-0.10$ for $\bsx \in D$, where $D$ is the disc of radius $1/4$ centered at the origin. A significant improvement in accuracy as well as run time is observed. For instance, the data in Table~\ref{AH} for $\kappa=120 \pi$ show that AH takes $77$ seconds to produce an approximate solution with three-digits of accuracy while FSPT achieves  better accuracy in just 2 seconds.

To make a similar comparison of with BH, we take $u^{i}(\bsx)=\exp (2\pi i x_1)$ and $m(\bsx)=-1$ for $\bsx \in D$, where $D$ is the unit disc centered at the origin. 
The corresponding numerical results for BH are taken from Table 5.2 on page 64 in \cite{hyde2002fast}. In the comparative study reported in Table~\ref{BH}, while the errors show a similar accuracy level for both schemes, a significant improvement in the run time is observed for the Fourier smoothed scheme. For instance, to achieve six-digits accuracy, BH requires $808$ seconds while the proposed FSPT scheme produces a similar level of accuracy in just $15$ seconds.

These two comparisons clearly show that our approach competes very well with some of the other existing alternatives.  

\begin{exmp}(\textit{Scattering by a non-smooth scattering geometries})
\end{exmp}

\begin{figure}[t!] 
	\begin{center}		
	\begin{minipage}[b]{0.315\linewidth}
	\centering
	\includegraphics[width=\linewidth]{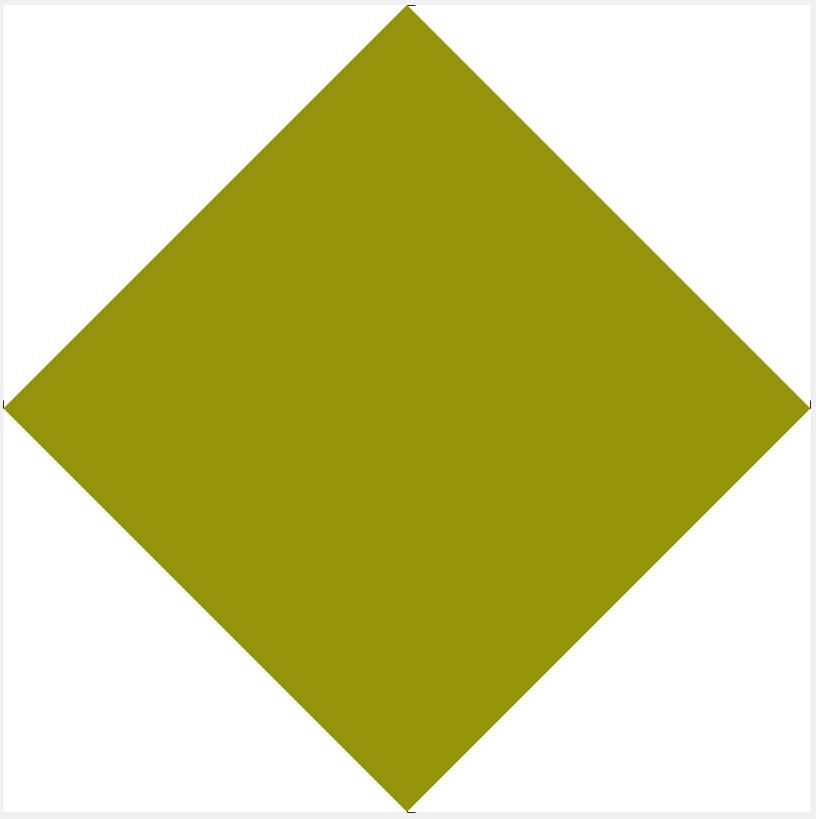}
	\subcaption{A scatterer with corners}
	\end{minipage}
	\hspace{0.05cm}	
	\begin{minipage}[b]{0.32\linewidth}
	\centering
	\includegraphics[width=\linewidth]{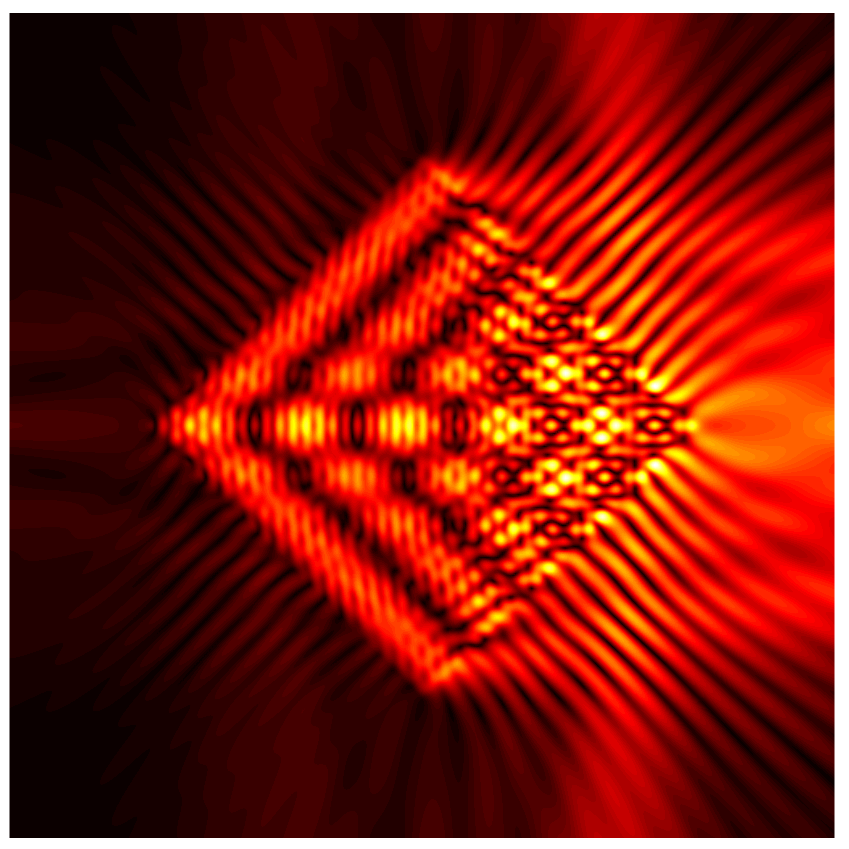}
	\subcaption{$|{u}^s|$}
	\end{minipage}
		\hspace{0.05cm}
		\begin{minipage}[b]{0.323\linewidth}
		\centering
		\includegraphics[width=\linewidth]{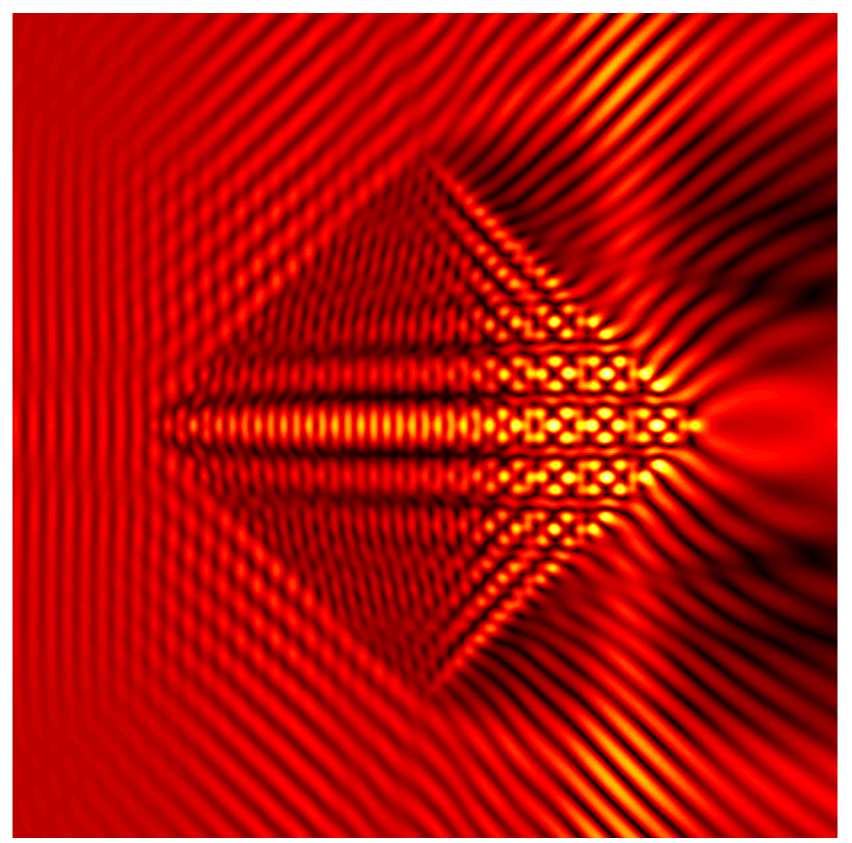}
		\subcaption{  $|u|$}
	\end{minipage}
	
	\caption{Scattering of a plane wave $ \exp (i \kappa x_1)$ by a square-shaped penetrable scatterer of acoustical size $\kappa \mathrm{d}=50$ and
	$m(\bsx) =-.90$ when $\bsx \in D$. The results correspond to a computational grid of size $512 \times 512 $ where a four-digits accuracy is achieved in the $\epsilon_{\infty}$ measure.}	
	\label{fig:-Square-K50}
	\end{center}
	\end{figure}

\begin{figure}[h!] 
	\begin{center}	
	\begin{minipage}[b]{0.315\linewidth}
		\centering
		\includegraphics[width=\linewidth]{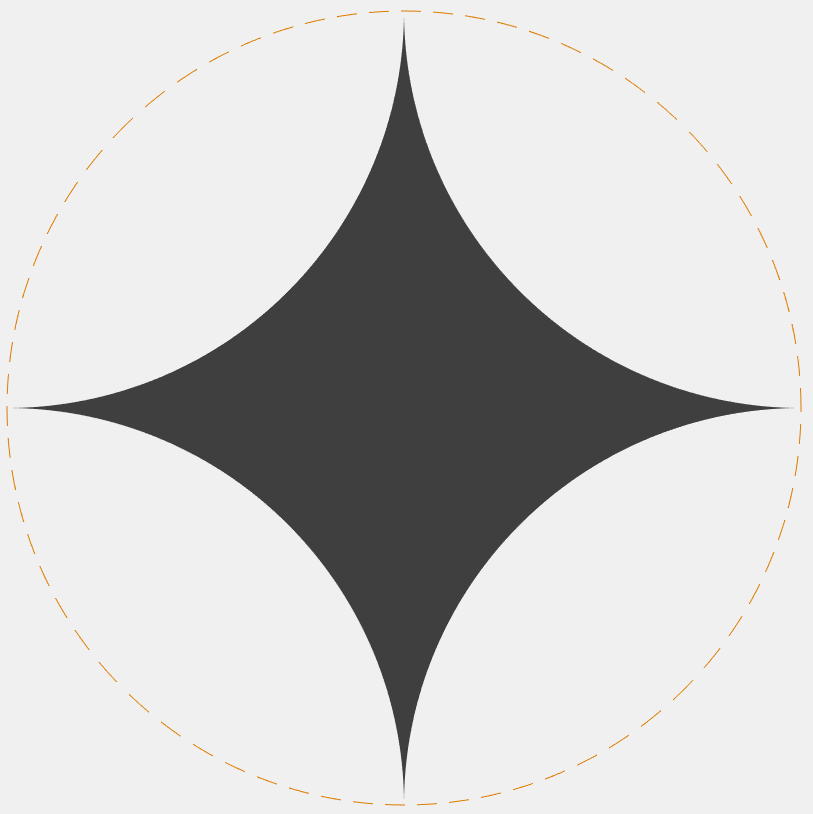}
		\subcaption{A scatterer with cusps}
	\end{minipage}
	\hspace{0.05cm}
	\begin{minipage}[b]{0.315\linewidth}
		\centering
		\includegraphics[width=\linewidth]{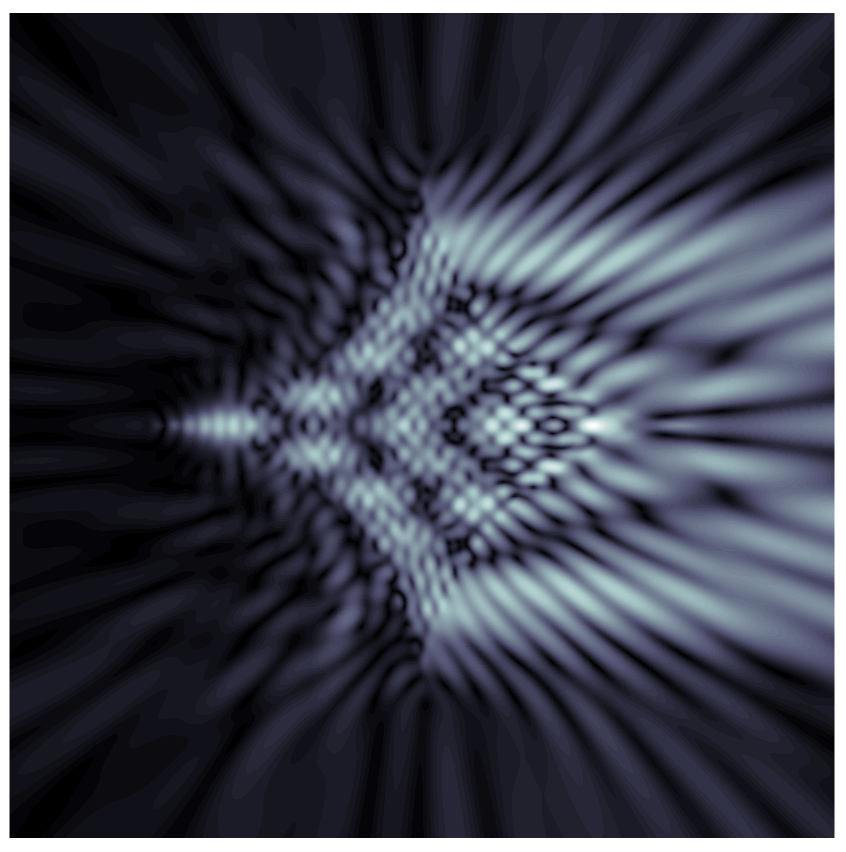}
		\subcaption{$|{u}^s|$}
	\end{minipage}
	\hspace{0.05cm}
		\begin{minipage}[b]{0.315\linewidth}
		\centering
		\includegraphics[width=\linewidth]{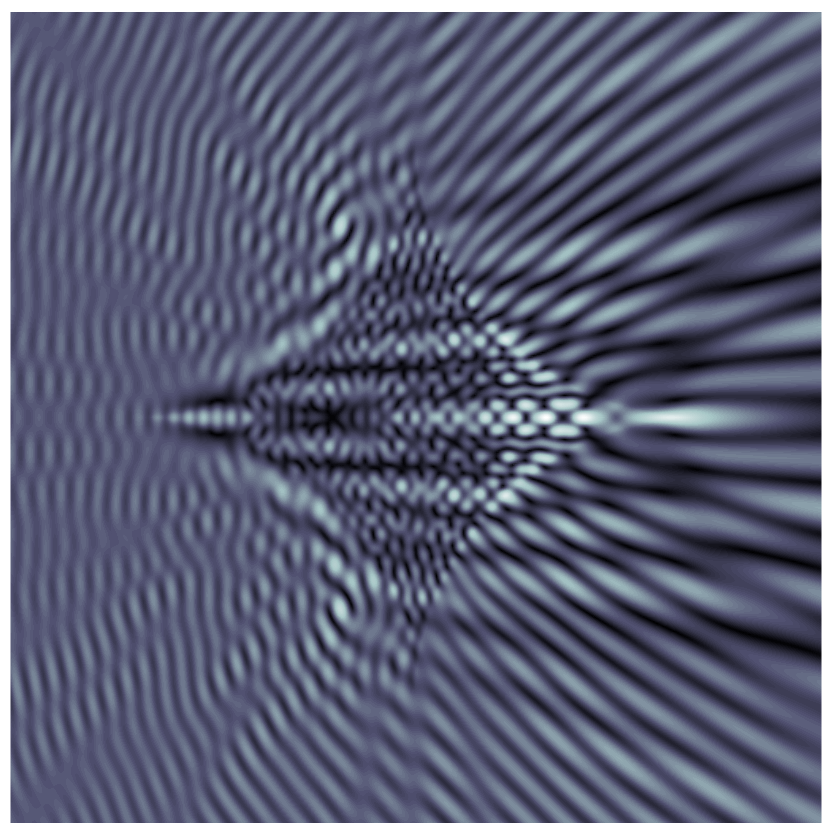}
		\subcaption{$|u|$}
	\end{minipage}
		
		\caption{Scattering of a plane wave $ \exp (i \kappa x_1)$ by a penetrable scatterer with cusps of acoustical size $\kappa \mathrm{d} =100$ and $m(\bsx) =-1.0$ when $\bsx \in D$. A computational grid of size $512 \times 512 $ is used that corresponds to three-digits accuracy in the computed solution.}
		\label{fig:-Cusp-K50}	
	\end{center}
\end{figure}

To show that the proposed methodology is not restricted only to problems with simple scattering geometries or with constant material properties, we next consider two additional  geometries, namely, a  square-shaped scatterer with corner singularities (see Figure~\ref{fig:-Square-K50}a) and a star-shaped scatterer containing cusps (shown in Figure~\ref{fig:-Cusp-K50}a). A plane wave incidence $u^{i}(\bsx)=\exp(i\kappa x_{1})$ is used for the scattering simulations. 

\begin{table}[b!]
	\begin{center} 
		\begin{tabular}{c|c|c|c|c|c|c|c} 
		\hline
		\hline		
		$n$& \multicolumn{2}{c|}{pre-corrected} &\multicolumn{5}{c}{Fourier smoothed pre-corrected}  \\ 
		\cline{2-8} &$\varepsilon_{2} $& $\varepsilon_{\infty} $&F &  $\varepsilon_{2}$& noc & $\varepsilon_{\infty}$& noc \\ 
		\hline
		\hline	
		
		$16$&$4.6\times 10^{-1}$&$8.6 \times 10^{-1}$&8&$4.6\times 10^{-1}$&-& $1.0\times 10^{0}$& -\\ 	
		
	   $32$&$4.4\times 10^{-1}$&$3.8 \times 10^{-1}$&16 &$4.6 \times 10^{-1}$&3.8&$8.4\times 10^{-1}$& 0.8\\ 

			$64$& $2.1\times 10^{-2}$&$5.1\times 10^{-2}$& 32 &$3.7 \times 10^{-3}$& 6.9&$1.7 \times 10^{-2}$&5.7 \\ 

			$128$ & $1.0\times 10^{-2}$&$4.5\times 10^{-2} $&64&$1.3 \times 10^{-4}$& 4.8&$1.2 \times 10^{-3}$& 3.7\\ 
		
		$256$&$7.5\times 10^{-3}$&$4.3\times 10^{-2} $&128  &$7.6 \times 10^{-6}$&4.1 &$2.2 \times 10^{-4}$&2.5\\

			$512$ &$2.6\times 10^{-3}$ &$1.5 \times 10^{-2}$& 256 &$3.0 \times 10^{-6}$ &1.3&$3.3\times 10^{-5} $&2.7\\ 
			
			$1024$ & $2.4\times 10^{-3}$&  $5.2 \times 10^{-3}$  &512& $4.5 \times 10^{-7}$&2.8&$7.0 \times 10^{-6}$&2.2\\ 	
			
			\hline	
			\hline
		\end{tabular}  
		\caption{Error report for the scattering by geometry with corners; reported results are computed for $u^{i}(\bsx)=\exp(i\kappa x_{1})$, $\kappa \mathrm{d}=50$ and
			$m(\bsx)=(x_{1}^2+x_{2}^2)\exp(-x_{1}^2-x_2^2)$ when $\bsx$ lies in the inhomogeneity.} 
		\label{tabel:-Corner} 
	\end{center} 
\end{table} 

\begin{table}[t!]
	\begin{center} 
		\begin{tabular}{c|c|c|c|c|c|c|c} \hline
			\hline
			$n$ & \multicolumn{2}{c|}{pre-corrected} &\multicolumn{5}{c}{Fourier smoothed pre-corrected}  \\ 
		\cline{2-8} &$\varepsilon_{2} $& $\varepsilon_{\infty} $&F &  $\varepsilon_{2}$& noc & $\varepsilon_{\infty}$& noc \\ 
			\hline
			\hline
			
			$16$& $5.4 \times 10^{-1} $& $1.1 \times 10^{0}$& 8 &$5.2 \times 10^{-1}$ &-&$1.1 \times 10^{0}$& -\\ 
			$32$&$5.2 \times 10^{-1}$& $1.0 \times 10^{0}$  &16 &$4.9 \times 10^{-1}$&-&$1.0 \times 10^{0}$& -\\ 
			$64$ &$2.2 \times 10^{-1}$  & $6.2 \times 10^{-1}$ &32  &$3.4 \times 10^{-2}$&3.9&$9.6\times 10^{-2}$&3.4\\ 
			$128$ &$7.8 \times 10^{-2}$&$1.9 \times 10^{-1}$&64   &$6.3 \times 10^{-4}$&5.7&$5.1\times 10^{-3}$& 4.2\\ 
			
			$256$ &$2.7\times 10^{-2}$&$7.5\times 10^{-2}$ &128 &$6.8\times 10^{-5}$&3.2&$1.0\times 10^{-3}$&2.3\\ 
			
			$512$  &$1.1\times 10^{-2}$  & $3.3\times 10^{-2}$ & 256&$1.2\times 10^{-5}$&2.5&$2.6\times 10^{-4}$&2.0\\ 
			
			$1024$  &$3.5\times 10^{-3}$&  $1.1\times 10^{-2}$ & 512&$2.2\times 10^{-6}$&2.5&$6.2\times 10^{-5}$&2.1\\ 				
			\hline	
			\hline
		\end{tabular}  
		\caption{Error report for the scattering by a geometry with cusps.  The incident field $u^{i}(\bsx)=\exp(i\kappa x_1)$, $\kappa \mathrm{d}=50$ and $m(\bsx)=-0.4$ when $\bsx \in D$ is used for this simulation.} 
		\label{tabel:-Cusp} 
	\end{center} 
\end{table} 

In Table~\ref{tabel:-Corner}, we present computational results for the scatterer depicted in Figure~\ref{fig:-Square-K50}a corresponding to the acoustic size $\kappa \mathrm{d}=50$ (the diameter of inhomogeneity $D$ is $\mathrm{d}$) and the contrast 
function within the scatterer given by the expression $m({\bf x})=(x_{1}^2+x_{2}^2)\exp(-x_{1}^2-x_{2}^2)$. For a pictorial illustration, absolute values of a four-digits accurate $u$ and $u^s$ are displayed in Figure~\ref{fig:-Square-K50}. On the other hand, Table~\ref{tabel:-Cusp} presents numerical results for the star-shaped scatterer displayed in Figure~\ref{fig:-Cusp-K50}a which essentially is the region enclosed by the four unit discs centered at  $(1,1), (1,-1), (-1,1)$ and $(-1,-1)$. The data in Table~\ref{tabel:-Cusp} correspond to the numerical simulation with $\kappa \mathrm{d}=50$ and $m({\bf x})=-0.4$ within the scatterer. A similar simulation with $m = -1$ within the scatterer is shown in Figure~\ref{fig:-Cusp-K50} where the depicted scattered field $u^{s}$ and the total field $u$ are computed to three-digits of accuracy.

\section{Conclusion}	\label{CL}
In this short communication, we have proposed a Fourier smoothed pre-corrected trapezoidal method that produces second order convergent solutions of the Lippmann-Schwinger integral equation for problems with a discontinuous refractive index at the interface and that has a computational complexity of $O(N \log N)$. The overall numerical scheme is very easy to implement and is applicable to a wide range of scattering geometries including those with singular geometric features such as corners and cusps. While we presented the method in the context of enhancing the rate of convergence of the pre-corrected trapezoidal rule, the proposed Fourier smoothing idea can be easily incorporated into other fast integral equation based solvers to enhance their rate of convergence when simulating scattering with discontinuous material properties. 

\section*{Acknowledgments}
AA acknowledges the support by Science \& Engineering Research Board through File No MTR/2017/000643.

\bibliographystyle{plain} 
\bibliography{References_WFS_LS_Solver}
	
\end{document}